\documentclass[12pt,reqno]{amsart}

\usepackage{amscd}

\newcommand{\N}{{\mathbb N}}
\newcommand{\Z}{{\mathbb Z}}
\newcommand{\Q}{{\mathbb Q}}
\newcommand{\C}{{\mathbb C}}
\newcommand{\R}{{\mathbb R}}
\renewcommand{\P}{{\mathbb P}}

\newcommand{\HH}{{\mathcal H}}

\newcommand{\NN}{{\mathcal N}}

\newcommand{\OO}{{\mathcal O}}

\newcommand{\ddd}{{\rm d}}
\newcommand{\www}{\widetilde}
\newcommand{\oooo}{\overline}

\renewcommand{\Im}{{\rm Im}}

\newcommand{\paa}{\partial}
\newcommand{\dtu}{{\partial_t^{-1}}}
\newcommand{\dt}{\partial_t}

\newcommand{\nnn}{\nabla}
\newcommand{\Var}{{\it Var}}
\newcommand{\Can}{{\it Can}}

\DeclareMathOperator{\Hom}{Hom}
\DeclareMathOperator{\id}{id}

\DeclareMathOperator{\Rad}{Rad}

\DeclareMathOperator{\Res}{Res}
\DeclareMathOperator{\res}{res}

\theoremstyle{plain}

\theoremstyle{definition}

\begin{document}

\title[Singularities and bilinear forms: a survey]
{Bilinear forms and hermitian forms for singularities:
a survey} 

\author{Claus Hertling }


\address{Claus Hertling\\
Lehrstuhl f\"ur algebraische Geometrie\\
Universit\"at Mannheim\\
B6, 26\\
68159 Mannheim\\
Germany}

\email{hertling\char64 math.uni-mannheim.de}

\maketitle

\section{Introduction}\label{s1}
\setcounter{equation}{0}

\noindent 
This is an english translation of the french article
\cite{He3}. It is made only for the arXiv. Only the
references \cite{BH}, \cite{He3} and \cite{De} are new. 

A one day conference on bilinear forms and hermitian forms
around isolated hypersurface singularities was organized
in Nancy on November 18, 2003. 
This article arose from the talk which I gave at the 
conference. Its aim is to collect different approaches
and view points and to give the relations between
them in precise formulas.
It does not contain new results.

The following diagram shows the different forms
which are considered in this article and indicates
their connections.

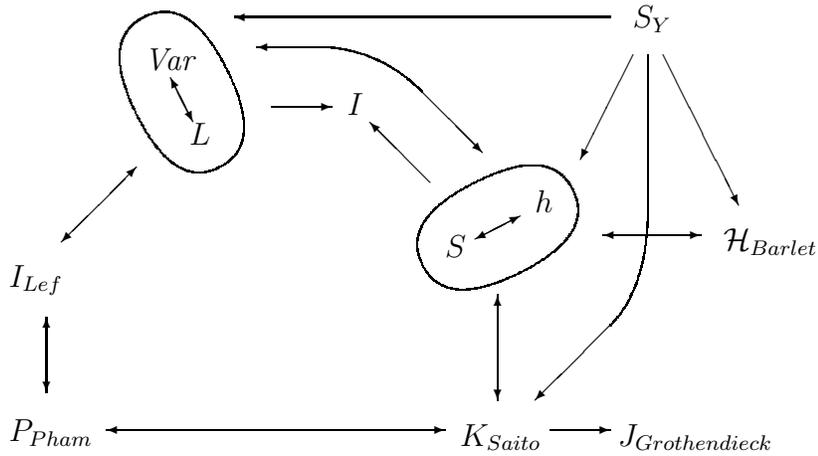
\setlength{\unitlength}{1cm}
\noindent
\begin{picture}(12,8)
\put(0.5,1){\makebox(0,0)[bl]{$P_{Pham}$}}
\put(1,2.2){\vector(0,1){0.5}}
\put(1,2.2){\vector(0,-1){0.5}}
\put(0.5,3){\makebox(0,0)[bl]{$I_{Lef}$}}
\put(1.7,4.2){\vector(-1,-1){0.5}}
\put(1.7,4.2){\vector(1,1){0.5}}
\put(2.3,6){\makebox(0,0)[bl]{\Var}}
\put(2.9,5){\makebox(0,0)[bl]{$L$}}
\put(2.75,5.7){\vector(1,-2){0.2}}
\put(2.85,5.5){\vector(-1,2){0.2}}
\bezier{100}(2.2,5.4)(2.7,4.4)(3.3,4.7)
\bezier{100}(3.4,6.0)(3.9,5.0)(3.3,4.7)
\bezier{100}(2.2,5.4)(1.7,6.4)(2.3,6.7)
\bezier{100}(3.4,6.0)(2.9,7.0)(2.3,6.7)

\put(4,5.5){\vector(1,0){0.8}}
\put(5,5.4){\makebox(0,0)[bl]{$I$}}
\put(6.1,4.5){\vector(-1,1){0.8}}
\put(6.3,3.5){\makebox(0,0)[bl]{$S$}}
\put(7.5,4.1){\makebox(0,0)[bl]{$h$}}
\put(6.9,3.85){\vector(2,1){0.4}}
\put(7.1,3.95){\vector(-2,-1){0.4}}
\bezier{100}(6.7,4.5)(5.7,4.0)(6.0,3.4)
\bezier{100}(6.7,4.5)(7.7,5.0)(8.0,4.4)
\bezier{100}(7.3,3.3)(8.3,3.8)(8.0,4.4)
\bezier{100}(7.3,3.3)(6.3,2.8)(6.0,3.4)

\bezier{100}(4.6,6.3)(5.4,6.3)(6.0,5.7)
\put(4.6,6.3){\vector(-1,0){0.8}}
\put(6.0,5.7){\vector(1,-1){0.8}}

\put(8.5,6.7){\vector(-1,0){5}}
\put(8.8,6.5){\makebox(0,0)[bl]{$S_Y$}}
\put(8.8,6.2){\vector(-1,-2){0.7}}

\put(9.2,6.2){\vector(1,-2){1}}
\put(10.0,3.6){\makebox{$\HH_{Barlet}$}}
\put(9.0,3.8){\vector(1,0){0.7}}
\put(9.0,3.8){\vector(-1,0){0.6}}

\put(9.0,6.2){\line(0,-1){2}}
\bezier{100}(9.0,4.2)(9.0,3.1)(8.5,2.6)
\put(8.5,2.6){\vector(-1,-1){1.0}}

\put(7.0,2.5){\vector(0,1){0.5}}
\put(7.0,2.5){\vector(0,-1){0.9}}

\put(3,1.2){\vector(1,0){3.3}}
\put(3,1.2){\vector(-1,0){1.2}}
\put(6.5,1.0){\makebox{$K_{Saito}$}}
\put(7.7,1.2){\vector(1,0){0.8}}
\put(8.6,1.0){\makebox{$J_{Grothendieck}$}}

\end{picture}

\noindent Here
\begin{enumerate}
\item $I$ is the intersection form on the homology of the 
Milnor fiber;
\item $L$ is the Seifert form;
\item $\Var$ is the variation operator;
\item $I_{Lef}$ is the intersection form for
Lefschetz thimbles;
\item $P_{Pham}$ is an induced form on the cohomologyy
bundle of Lefschetz thimbles; it was defined by F. Pham;
\item $S$ is a polarizing form for Steenbrink's mixed
Hodge structure;
it is defined on the cohomology of a Milnor fiber;
\item $h$ is a hermitian form on the same space;
it was defined by D. Barlet;
\item $S_Y$ is the polarizing form on the cohomology
of a fiber of a compactification of the Milnor fibration;
\item $\HH_{Barlet}$ is a hermitian form on the
Brieskorn lattices, defined by D. Barlet;
\item $K_{Saito}$ 
is the form of K. Saito on the Brieskorn lattice;
\item and finally, $J_{Grothendieck}$ 
is a form which is induced by the Grothendieck residue.
\end{enumerate}

\noindent 
Except for $\Var$, these are bilinear or hermitian forms.\\
In the diagram, the monodromy $M$ is missing.
The Seifert form $L$ and the variation operator $\Var$
determine the monodromy.\\
But the monodromy must be taken for granted in the 
correspondences which are indicated by the arrows
in the diagram.

\noindent Chapter \ref{s2} describes the forms of
topological origin and the relations between them.
Chapter \ref{s3} describes the transcendent forms
and the additional relations.
In this paper $\N=\{0,1,2,...\}$. 

\noindent
I thank Daniel Barlet.

\section{Topological forms}\label{s2}
\setcounter{equation}{0}

\subsection{Milnor fibration and monodromy}\label{s2.1}

One fixes a germ
\\  $f:(\C^{n+1},0)\to (\C,0)$ 
of a holomorphic function with an isolated singularity 
at 0 which has Milnor number $\mu$. 
A standard reference for the topological properties of the
Milnor fibration, which are recalled in the sections
\ref{s2.1} to \ref{s2.4} is 
\cite[Part I]{AGV}.

\noindent 
One obtains a representative by choosing 
 $\varepsilon$ and $\eta$  with $0<\eta\ll\varepsilon\ll 1$
and by defining $\Delta=\Delta_\eta=\{z\in \C\ |\ |z|<\eta\}$,
$B_\varepsilon=\{x\in \C^{n+1}\ |\ \|x\| <\varepsilon\}$, and 
$X=f^{-1}(\Delta)\cap B_\varepsilon$.

\noindent 
The fibers $X_t=f^{-1}(t)\subset X$, 
$t\in \Delta^*=\Delta-\{0\}$, 
are smooth and have the same homotopy type as a bouquet
of $\mu$ spheres of dimension $n$. The fiber $X_0$ is
contractible.

\noindent 
The Milnor fibration $f:X-X_0\to \Delta^*$ 
is locally trivial (as a $C^{\infty}$ map).
One has a monodromy on the homology $H_n(X_t,\Z)\cong \Z^\mu$ 
($t\neq 0$) and also on $H^n(X_t,\Z)$, 
which we denote always by $M$.

It is quasiunipotent (see \cite{Br}).
It decomposes on $H_n(X_t,\Q)$ into a semisimple part
$M_s$ and a unipotent part $M_u$, and one has
$M=M_s\cdot M_u=M_u\cdot M_s$.
The nilpotent part $N:=\log M_u$ acts on 
$H_n(X_t,\Q)$ and satisfies $N^{n+1}=0$. 
We will use the notations
\begin{eqnarray}
H_n(X_t,\C)_\lambda :=\ker (M_s-\lambda\cdot \id)\subset H_n(X_t,\C),
\label{2.1}\\
H_n(X_t,\Q)_{\neq 1}:=H_n(X_t,\Q)\cap 
\bigoplus_{\lambda\neq 1}H_n(X_t,\C)_\lambda ,\label{2.2}.
\end{eqnarray}
We have with the same notation for $H_n(X_t,\Q)_1$
\begin{eqnarray} 
H_n(X_t,\Q)=H_n(X_t,\Q)_{\neq 1}\oplus H_n(X_t,\Q)_1.\label{2.3}
\end{eqnarray}
We will use analogous notations for the cohomology.

\subsection{Intersection form $I$}\label{s2.2}

The intersection form 
\begin{eqnarray}
I:H_n(X_t,\Z)\times H_n(X_t,\Z)\to \Z \label{2.4}
\end{eqnarray}
is $(-1)^n$-symmetric and $M$-invariant. Its radical is 
\begin{eqnarray}
\Rad I=\ker (M-\id)\subset H_n(X_t,\Z).\label{2.5}
\end{eqnarray}
This equality is a consequence of the long exact sequence
for the pair
$(X_t,\paa X_t)$:
\begin{eqnarray}
0&\to& H_n(\paa X_t,\Z) \to H_n(X_t,\Z) 
\stackrel{\Can}{\longrightarrow} H_n(X_t,\paa X_t,\Z) \label{2.6}\\
&\to& H_{n-1}(\paa X_t,\Z) \to 0\nonumber
\end{eqnarray}
and of the isomorphism
$H_n(\paa X_t,\Z)\cong \ker(M-\id)\subset H_n(X_t,\Z)$.
The intersection form and the morphism $\Can$ 
induce one another by the formula
\begin{eqnarray}
I(a,b)=\langle \Can(a),b\rangle.\label{2.7}
\end{eqnarray}
In this paper, $\langle ,\rangle$ denotes the canonical
bilinear form between two dual spaces;
in the moment, they are 
$H_n(X_t,\paa X_t,\Z)$ and $H_n(X_t,\Z)$.

\subsection{Variation operator $\Var$}\label{s2.3}
The intersection form between (absolute) cycles and
relative cycles is a duality between the spaces
$H_n(X_t,\Z)$ and $H_n(X_t,\paa X_t,\Z)$.
We have therefore a canonical isomorphism
\begin{eqnarray}
H^n(X_t,\Z)\cong H_n(X_t,\paa X_t,\Z).\label{2.8}
\end{eqnarray}
The variation operator is the morphism
\begin{eqnarray}
\Var: H^n(X_t,\Z)\to H_n(X_t,\Z)\label{2.9}
\end{eqnarray}
which is defined by
$\Var([\gamma]):= [Mon(\gamma)-\gamma]$,
where $\gamma$ 
is a representative of the cycle 
$[\gamma]\in H_n(X_t,\paa X_t,\Z)$,
and where $Mon$ 
denotes a representative of the monodromy on the fiber$X_t$,
which acts trivially on $\paa X_t$. 
Therefore $Mon(\gamma)-\gamma$ 
represents a cycle without boundary.

The variation operator $\Var$ is an isomorphism.
A simple way to see this consists in comparison with the
intersection form  $I_{Lef}$ of Lefschetz thimbles,
see \ref{s2.5}. By definition of $\Var$, 
the following diagram commutes.

\begin{eqnarray}\label{2.9b}
\textup{
\begin{picture}(11,4)
\put(3,1){\makebox(0,0){$H_n(X_t,\Z)$}}
\put(9,1){\makebox(0,0){$H^n(X_t,\Z)$}}
\put(3,3){\makebox(0,0){$H_n(X_t,\Z)$}}
\put(9,3){\makebox(0,0){$H^n(X_t,\Z)$}}
\put(3,2.5){\vector(0,-1){1}}
\put(2.9,2){\makebox(0,0)[r]{$M-\id$}}
\put(9,2.5){\vector(0,-1){1}}
\put(9.1,2){\makebox(0,0)[l]{$M-\id$}}
\put(4.2,1){\vector(1,0){3.6}}
\put(6,1.3){\makebox(0,0){\Can}}
\put(4.2,3){\vector(1,0){3.6}}
\put(6,3.3){\makebox(0,0){\Can}}
\put(7.8,2.7){\vector(-3,-1){3.6}}
\put(6,2.4){\makebox(0,0){$\Var$}}
\end{picture}
}
\end{eqnarray}

The fact that $\Var$ is an isomorphism shows \eqref{2.5},
\begin{eqnarray*}
\ker (M-\id)= \ker \Can =\Rad I.
\end{eqnarray*}
On $H^n(X_t,\Q)$ we have
\begin{eqnarray}
\Im (M-\id)&=&\Im (\Can)\label{2.10}\\
&=&H^n(X_t,\Q)_{\neq 1}\oplus (H^n(X_t,\Q)_1\cap\Im (\Can)). \nonumber
\end{eqnarray}
We can define on this space a form $I^{coh}$ 
by the formula
\begin{eqnarray}
I^{coh}(A,B):= I(\Can^{-1}A,\Can^{-1}B). \label{2.11}
\end{eqnarray}
In fact, we do not have unique preimages of $A$ and $B$
by $\Can$. But the different preimages differ only by
elements of $\Rad I$. 
In the formula \eqref{2.12}, we have a similar inessential
ambiguity. On $H^n(X_t,\Q)_{\neq 1}$ we have no ambiguity; 
$\Can$ and $M-\id$ are invertible.

The form $I^{coh}$ is $(-1)^n$-symmetric and  $M$-invariant; 
it is nondegenerate on $H^n(X_t,\Q)_{\neq 1}$. 
The diagram shows that we have for  $A,B\in \Im(\Can)$
\begin{eqnarray}
I^{coh}(A,B):= \langle A,\Can^{-1}B\rangle 
= \langle A,\Var\circ \frac{1}{M-\id}(B)\rangle . \label{2.12}
\end{eqnarray}

\subsection{Seifert form $L$}\label{s2.4}
The Seifert form is a form 
\begin{eqnarray}
L:H_n(X_t,\Z)\times H_n(X_t,\Z)\to \Z,\label{2.13}
\end{eqnarray}
which is unimodular. it is neither symmetric nor 
antisymmetric; it satisfies the following properties, 
\begin{eqnarray}
L(M\, a,b)&=&(-1)^{n+1}L(b,a), \label{2.14}\\
I(a,b)&=& L((M-\id)a,b) \label{2.15}\\
&=& -L(a,b) + (-1)^{n+1}L(b,a) \nonumber.
\end{eqnarray}
The Seifert form and the variation operator 
determine one another by the formula
\begin{eqnarray}
L(a,b)= \langle \Var^{-1}\, a,b\rangle. \label{2.16}
\end{eqnarray}
Therefore one can read the properties \eqref{2.14} and 
\eqref{2.15} as properties of the variation operator.
A simple way to understand \eqref{2.14} and 
\eqref{2.16} is by comparison with the intersection form
$I_{Lef}$ of Lefschetz thimbles, see \ref{s2.5}.
The relation \eqref{2.15} is a consequence of \eqref{2.16}, 
$\Can=\Var^{-1}\circ(M-\id)$ and \eqref{2.7}. 

In the literature, there are different conventions for
the definition of the Seifert form. Here we follow
\cite[Part I]{AGV}, where one finds also the proofs of 
\eqref{2.14} -- \eqref{2.16}.

For the definition, we need several notations.
Let $\pi:\C\to\C^*$ be the universal covering,
where $\pi(\zeta)=e^{2i\pi \zeta}$, and let  
$\pi:\www{\Delta^*}\to \Delta^*$ be its restriction. 
The choice of a trivialization of the fibration
$\pi^* (X-X_0)\to \www{\Delta^*}$ defines a diffeomorphism
\begin{eqnarray}
Mon_{[\beta]}&:&X_t\to X_{t \cdot\pi(\beta)}, 
\quad \beta\in \R\label{2.17}\\
&& X_t \cong (\pi^* X)_\theta \cong (\pi^*X)_{\theta+\beta} 
\cong X_{t\cdot\pi(\beta)}, \quad \hbox{ with }
\pi(\theta)=t, \nonumber
\end{eqnarray}
and an induced isomorphism
\begin{eqnarray}
M_{[\beta]}:H_n(X_t,\Z)\to H_n(X_{t \cdot\pi(\beta)},\Z).\label{2.18}
\end{eqnarray}
In particular, $M_{[k]}=M^k$ if $k\in\Z$.
Following Milnor, one can also choose the diffeomorphisms 
\begin{eqnarray}
\widetilde{.}:X_t\to \www {X_t}:=(f/|f|)^{-1}(t/|t|) \subset S^{2n+1}-X_0.
\label{2.19}
\end{eqnarray}
For the $n$-cycles of $S^{2n+1}$, one has a linking form.
One defines the Seifert form by 
\begin{eqnarray}
L([a],[b]):= \hbox{ (linking form) }
(\www a, \www{Mon_{[1/2]}b}).\label{2.20}
\end{eqnarray}
One can find an oriented $(n+1)$-chain $A$, 
with boundary $\paa A=\www a$ in $S^{2n+1}$
such that the orientation of $A$ 
is given near $\paa A$ by
the pair $(v,\textup{orientation of }\www a)$, where $v$ 
is a vector transverse to $\paa A$ 
which points to the exterior. 
The linking form is defined in a way such that
\begin{eqnarray}
L([a],[b])= \hbox{ (intersection number) }(A,\www{Mon_{[1/2]}b}) \label{2.21}
\end{eqnarray}

\subsection{Lefschetz thimbles, the form $I_{Lef}$}\label{s2.5}
The Milnor fibration extends to $f:\oooo X\to \oooo \Delta$ 
with the fibers $\oooo X_t$, $t\in \oooo\Delta$; 
here $\oooo X$ and $\oooo\Delta$
are the closures of $X$ in $\C^{n+1}$ and of 
$\Delta=\Delta_\eta$ in $\C$.
The long exact sequence of the pair 
$(\oooo X,\oooo X_{\eta\cdot\theta})$
with $\theta\in S^1$ is
\begin{eqnarray} \label{2.22}
0=H_{n+1}(\oooo X,\Z)&\to& H_{n+1}(\oooo X,\oooo X_{\eta\theta},\Z) \\
&\stackrel{\paa_*}{\longrightarrow}& H_n(\oooo X_{\eta\theta},\Z)
\to H_n(\oooo X,\Z)=0 .\nonumber
\end{eqnarray}
Therefore $\paa_*$ is a canonical isomorphism. \\
The elements of $ H_{n+1}(\oooo X,\oooo X_{\eta\theta},\Z) $
are represented by the Lefschetz thimbles (e.g. \cite{Ph1}).
A Lefschetz thimble $\Gamma$ is a family
$\bigcup_{t\in \Im{(\sigma)}}\delta(t)$ of $n$-cycles 
$\delta(t)\subset X_t$
over a path $\sigma:[0,1]\to \oooo\Delta$ 
which is a sufficiently smooth embedding with 
$\sigma(0)=0$ and $\sigma(1)=\eta\theta$.

One has an intersection form $I_{Lef}$ for the Lefschetz
thimbles with boundaries in $\oooo X_{\eta\theta}$ and 
$\oooo X_{-\eta\theta}$,
\begin{eqnarray} \label{2.23}
I_{Lef}:H_{n+1}(\oooo X,\oooo X_{\eta\theta},\Z)\times
H_{n+1}(\oooo X,\oooo X_{-\eta\theta},\Z)\to \Z .
\end{eqnarray}
It is $(-1)^{n+1}$-symmetric.
It is unimodular (e.g. \cite{Ph1}, \cite{AGV}). 
A proof uses a deformation of $f$ to a function with $\mu$
singularities of type $A_1$ which have pairwise different
critical values. 
A system of $\mu$ paths $\sigma_1,...,\sigma_\mu$ which join
the $\mu$ critical values with $\eta\theta$ and which
do not intersect one another, give a distinguished basis
(modulo the signs) of Lefschetz thimbles for 
$H_{n+1}(\oooo X,\oooo X_{\eta\theta},\Z)$.
A complementary system of paths 
$\www \sigma_1,...\www \sigma_\mu$ which join the $\mu$
critical values with $-\eta\theta$, which do not intersect
one another and which do not meet $\sigma_1,...,\sigma_\mu$, 
give a dual basis (modulo the signs) of 
$H_{n+1}(\oooo X,\oooo X_{-\eta\theta},\Z)$ with respect
to $I_{Lef}$.
Of course, the topology of the pair 
$(\oooo X,\oooo X_{\eta\theta})$
does not change along the deformation of $f$.

With the convention of the orientation which was used in 
\ref{s2.4}, one sees without difficulty
(e.g. \cite[Part I]{AGV}) that 
\begin{eqnarray} \label{2.24}
L(a,b)= (-1)^{n+1}I_{Lef}
(\paa_*^{-1}(a),\paa_*^{-1}(M_{[1/2]}b)).
\end{eqnarray}
The $(-1)^{n+1}$-symmetry of $I_{Lef}$, \eqref{2.24}, 
the fact that $I_{Lef}$ is invariant with respect to 
$M_{[\beta]}$ ($\beta\in \R$), and
$M_{[1/2]}\circ M_{[1/2]}=M$ 
show immediately the formula \eqref{2.14}.

Two of the three formulas \eqref{2.16}, \eqref{2.24} and
\begin{eqnarray} \label{2.25}
I_{Lef}(\paa_*^{-1}(a),\paa_*^{-1}(M_{[1/2]}b)) 
= (-1)^{n+1}\langle \Var^{-1}(a),b\rangle
\end{eqnarray}
induce the third one.
One can prove \eqref{2.25} with the following observation. 
If $\alpha$ represents a relative cycle in  
$H_{n}(\oooo X_{\eta\theta},\paa \oooo X_{\eta\theta},\Z)$,
the union of $\bigcup_{\theta\in [0,1]}Mon_{[\theta]}(\alpha)$ 
with a suitable subset
 $\bigcup_{t\in \Delta}\paa X_t$ represents
the class $\paa^{-1}_*(\Var \alpha)$ in
$H_{n+1}(\oooo X,\oooo X_{\eta\theta},\Z)$. 
The sign $(-1)^{n+1}$ in \eqref{2.25} comes from the 
orientations.

\subsection{Cohomology bundle of the Lefschetz thimbles, 
$P_{Pham}$}\label{s2.6}
The groups 
$H_{n+1}(\oooo X,\oooo X_{\eta\theta},\Z)$, $\theta\in S^1$,
give a local system of $\Z$-modules which are free of rank 
$\mu$, on $S^1$,
which one can extend to $\C^*$ by the retraction 
$\C^*\to S^1,\ z\mapsto z/|z|$.
One obtains a dual local system 
\begin{eqnarray}\label{2.26}
H^{Lef}_\Z&=&\bigcup_{z\in \C^*}H^{Lef}_\Z(z)\quad \hbox{ avec }\\
H^{Lef}_\Z(z)&=& \Hom(H_{n+1}(\oooo X,\oooo X_{\eta z/|z|},\Z),\Z)
\nonumber
\end{eqnarray}
and a flat vector bundle 
$H^{Lef}_\C=\bigcup_{z\in \C^*}H^{Lef}_\C(z)$.
The unimodular form $I_{Lef}$ induces the isomorphisms
\begin{eqnarray}\label{2.27}
\pi_{Lef,z}: H_{n+1}(\oooo X,\oooo X_{\eta z/|z|},\Z)\to H^{Lef}_\Z(-z),
\end{eqnarray}
a unimodular form
\begin{eqnarray}\label{2.28}
&&I^{Lef}:H^{Lef}_\Z(z)\times H^{Lef}_Z(-z)\to \Z,\\
&&I^{Lef}=I_{Lef}\circ( \pi_{Lef,-z}^{-1},\pi_{Lef,z}^{-1})\nonumber
\end{eqnarray}
and a form
\begin{eqnarray}\label{2.29}
&&P_{Pham}:H^{Lef}_\Z(z)\times H^{Lef}_Z(-z)\to \frac{1}{(2i\pi)^{n+1}}\Z,\\
&&P_{Pham}= (-1)^{n(n+1)/2}\frac{1}{(2i\pi)^{n+1}}I^{Lef},\nonumber
\end{eqnarray}
which had been defined by Pham \cite{Ph2}.
It is flat on $H^{Lef}_\C$, $(-1)^{n+1}$-symmetric annd 
nondegenerate. After a Fourier-Laplace transformation,
it gives the form of K. Saito on the Brieskorn lattice,
see \ref{s3.8}.

\subsection{Compactification of the Milnor fibration, 
form $S_Y$}\label{s2.7}
Following \cite{Sch}, one can compactify the fibration 
$f:X\to \Delta$ to a fibration $f_Y:Y\to \Delta$ 
with the following properties
(modulo shrinking $\Delta$): 
\begin{enumerate}
\item the fibers $Y_t$ are hypersurfaces in $\P^{n+1}$, 
smooth for $t\neq 0$ and with a unique singularity in 
$0\in X_0\subset Y_0$ for $t=0$;
\item  one has an exact sequence 
\begin{eqnarray}\label{2.30}
0\to H^n(Y_0,\Q)\to H^n(Y_t,\Q) \stackrel{i^*}{\longrightarrow} 
H^n(X_t,\Q)\to 0, \quad t\neq 0.
\end{eqnarray}
\end{enumerate}

\noindent
One knows that $H^n(Y_0,\Q)\cong\ker(M_Y-\id)$, by the
theorem of invariant cycles.
Here $M_Y$ denotes the monodromy on the homology and the 
cohomology of the fibers of $f_Y$. 

\noindent
The intersection form $I_Y$ on $H_n(Y_t,\Q)$ is 
nondegenerate, $M_Y$-invariant and $(-1)^n$-symmetric.
It induces the isomorphisms
\begin{eqnarray}\label{2.30b}
\Can_Y:H_n(Y_t,\Q)\to H^n(Y_t,\Q) \textup{ (Poincar\'e duality)}
\end{eqnarray}
and a form
\begin{eqnarray}\label{2.30c}
I_Y^{coh}=I\circ (\Can_Y^{-1},\Can_Y^{-1})\textup{ on }
H^n(Y_t,\Q).
\end{eqnarray}

\noindent
The form
\begin{eqnarray}\label{2.31}
S_Y:= (-1)^{n(n-1)/2}I_Y^{coh}: H^n(Y_t,\Q)\times 
H^n(Y_t,\Q)\to \Q
\end{eqnarray}
is a polarizing form of pure Hodge structures on the groups
$H^n(Y_t,\C)$, $t\in\Delta^*$.

Following \cite{St}, one can enrich the exact sequence 
\eqref{2.30} with mixed Hodge structures (MHS),
which are compatible with the morphisms: the MHS of Deligne 
on $H^n(Y_0,\C)$, the limit MHS of Schmid on $H^n(Y_t,\C)$ 
and the MHS of Steenbrink on $H^n(X_t,\C)$. 
The MHS of Schmid on $H^n(Y_t,\C)$ is also polarized 
by the form $S_Y$, where the definition of a polarized 
MHS is analogous to the limit MHS of Schmid, see \cite{CK} 
or \cite{He1}.

\subsection{Polarizing form $S$}\label{s2.8}
The exact sequence \eqref{2.30} of MHS induces a polarizing
form $S$ on $H^n(X_t,\Q)$ for its MHS. 
In order to determine this form, one observes that the MHS
on $H^n(Y_t,\C)$ and $H^n(X_t,\C)$ are invariant by the 
semisimple parts of the monodromies $M_Y$ and $M$. 
Here one applies the notations of \ref{s2.1} also to the 
fibration $f_Y:Y\to \Delta$;
in particular, $N_Y=\log M_{Y,u}$ is the nilpotent part of the
monodromy $M_Y$. The exact sequence \eqref{2.30} gives an
isomorphism $ H^n(Y_t,\Q)_{\neq 1}\to H^n(X_t,\Q)_{\neq 1}$.
With the definition of $I^{coh}$ which was given in 
\eqref{2.11}, one obtains
\begin{eqnarray}\label{2.32}
S=(-1)^{n(n-1)/2}I^{coh} \cong (-1)^{n(n-1)/2}I_Y^{coh} 
\cong S_Y \  \hbox{ on }H^n(Y_t,\Q)_{\neq 1}.
\end{eqnarray}
The exact sequence \eqref{2.30} restricts to an exact sequence
\begin{eqnarray}\label{2.33}
0&\to& \ker(N_Y:H^n(Y_t,\Q)_1\to H^n(Y_t,\Q)_1)\\
&\to& H^n(Y_t,\Q)_1\to H^n(X_t,\Q)_1 \to 0.\nonumber
\end{eqnarray}
The definition of a polarized mixed Hodge structure
\cite{CK}\cite{He1} shows that for $a,b\in H^n(X_t,\Q)_1$ with
preimages $a_Y,b_Y\in H^n(Y_t,\Q)_1$ one has  
\begin{eqnarray}\label{2.34}
S(a,b)=S_Y(a_Y,(-N_Y)b_Y).
\end{eqnarray}
The right hand side is independent of the choices of 
$a_Y$ et $b_Y$; in fact, $N_Y$ is an infinitesimal isometry
of $S_Y$.
Therefore the polarizing form $S$ is $M$-invariant, 
nondegenerate, $(-1)^n$-symmetric on $H^n(X_t,\Q)_{\neq 1}$ and
$(-1)^{n+1}$-symmetric on $H^n(X_t,\Q)_1$.

One can also describe $S$ only in terms of the topology of the  
fibration $f:X\to \Delta$. For this purpose,
one can complete the diagram \eqref{2.9b} 
(with coefficients in $\Q$).

\begin{eqnarray}\label{2.34b}
\textup{
\begin{picture}(11,8)
\put(3,3){\makebox(0,0){$H_n(X_t,\Q)$}}
\put(9,3){\makebox(0,0){$H^n(X_t,\Q)$}}
\put(3,5){\makebox(0,0){$H_n(X_t,\Q)$}}
\put(9,5){\makebox(0,0){$H^n(X_t,\Q)$}}
\put(3,1){\makebox(0,0){$H_n(Y_t,\Q)$}}
\put(9,1){\makebox(0,0){$H^n(Y_t,\Q)$}}
\put(3,7){\makebox(0,0){$H_n(Y_t,\Q)$}}
\put(9,7){\makebox(0,0){$H^n(Y_t,\Q)$}}
\put(3,4.5){\vector(0,-1){1}}
\put(2.9,4){\makebox(0,0)[r]{$M-\id$}}
\put(9,4.5){\vector(0,-1){1}}
\put(9.1,4){\makebox(0,0)[l]{$M-\id$}}
\put(4.2,3){\vector(1,0){3.6}}
\put(6,3.5){\makebox(0,0){\Can}}
\put(4.2,5){\vector(1,0){3.6}}
\put(6,5.5){\makebox(0,0){\Can}}
\put(7.8,4.7){\vector(-3,-1){3.6}}
\put(6,4.4){\makebox(0,0){$\Var$}}
\put(3,2.5){\vector(0,-1){1}}
\put(2.9,2){\makebox(0,0)[r]{$i_*$}}
\put(9,1.5){\vector(0,1){1}}
\put(9.1,2){\makebox(0,0)[l]{$i^*$}}
\put(4.2,1){\vector(1,0){3.6}}
\put(6,1.5){\makebox(0,0){$\Can_Y$}}
\put(3,5.5){\vector(0,1){1}}
\put(2.9,6){\makebox(0,0)[r]{$i_*$}}
\put(9,6.5){\vector(0,-1){1}}
\put(9.1,6){\makebox(0,0)[l]{$i^*$}}
\put(4.2,7){\vector(1,0){3.6}}
\put(6,7.5){\makebox(0,0){$\Can_Y$}}
\end{picture}
}
\end{eqnarray}
One sees 
\begin{eqnarray}\label{2.35}
i_*\circ \Var \circ i^*\circ \Can_Y = M_Y-\id.
\end{eqnarray}
With this formula for the space $H^n(X_t,\Q)_1$ and with 
\eqref{2.12} for the space $H^n(X_t,\Q)_{\neq 1}$ 
one finds for $a,b\in H^n(X_t,\Q)$
\begin{eqnarray}\label{2.36}
S(a,b) = (-1)^{n(n-1)/2}\cdot \langle a, \Var\circ \nu(b)\rangle,
\end{eqnarray}
where the isomorphism $\nu: H^n(X_t,\Q)\to H^n(X_t,\Q)$ 
is $M$-invariant and 
\begin{eqnarray} \label{2.37}
\nu &=& \frac{1}{M-\id} \quad \hbox{ on } 
H^n(X_t,\Q)_{\neq 1},\\
\nu &=& \frac{-N}{M-\id} \quad \hbox{ on } H^n(X_t,\Q)_{1}.
\label{2.38}
\end{eqnarray}

Furthermore, the definition of $I^{coh}$ on 
$\Im (\Can)\cap H^n(X_t,\Q)_{1}$
shows that for $a,b\in\Im (\Can)\cap H^n(X_t,\Q)_{1}$
\begin{eqnarray}\label{2.39}
S(a,b)=(-1)^{n(n-1)/2}I^{coh}(a,(-N)b).
\end{eqnarray}

\subsection{Hermitian form $h$}\label{s2.9}
D. Barlet \cite{Ba} defined a sesquilinear form
($\C$-linear in the left argument, 
$\C$-semi-linear in the right argument), which is also 
hermitian and $M$-invariant,
\begin{eqnarray}\label{2.40}
h: H^n(X_t,\C)\times H^n(X_t,\C)\to \C.
\end{eqnarray}
F. Loeser \cite{Loe} observed that one can write this form
in terms of a compactification of the Milnor fibration,
and also in terms of the polarizing form. 
In fact, $S$ and $h$ determine one another via the following
formulas.
\begin{eqnarray}\label{2.41}
h(a,b)= (-1)^{n(n-1)/2}\frac{1}{(2i\pi)^n} S(a,\oooo b)
\quad \hbox { on } H^n(X_t,\C)_{\neq 1},\\
h(a,b)= (-1)^{n(n-1)/2}\frac{-1}{(2i\pi)^{n+1}} S(a,\oooo b)
\quad \hbox { on } H^n(X_t,\C)_{1}.\label{2.42}
\end{eqnarray}
In terms of a compactification of the Milnor fibration and 
the form $I_Y^{coh}$, the form $h$ is given by the formulas
\begin{eqnarray}\label{2.43}
h(a,b)= \frac{1}{(2i\pi)^n} I_Y^{coh}(a_Y,\oooo b_Y)
\quad \hbox { on } H^n(X_t,\C)_{\neq 1},\\
h(a,b)= \frac{1}{(2i\pi)^{n}} 
I_Y^{coh}((\frac{-N_Y}{2i\pi})a_Y,\oooo b_Y)
\quad \hbox { on } H^n(X_t,\C)_{1}\label{2.44};
\end{eqnarray}
here, $a_Y$ and $b_Y\in H^n(Y_t,\C)$ 
are preimages of $a$ and $b$.

The definition of Barlet uses integration over the fibers
of the Milnor fibration and the asymptotic expansion
of the sections of the Gauss-Manin connection, see \ref{3.4}.

\section{Transcendent forms}\label{s3}
\setcounter{equation}{0}

\subsection{Brieskorn lattice}\label{s3.1}
The following space had first been considered by Brieskorn 
\cite{Br},
\begin{eqnarray}\label{3.1}
E := \Omega^{n+1}_{\C^{n+1},0}/\ddd f \wedge \ddd \Omega^{n-1}_{\C^{n+1},0}. 
\end{eqnarray}
It is a free module of rank $\mu$ over the rings $\C\{t\}$ (
M. Sebastiani, B. Malgrange) and $\C\{\{\dtu\}\}$. 
One can identify $E$ with a subspace of the space of germs
at 0 of holomorphic sections of the cohomology bundle
\begin{eqnarray}\label{3.2}
H^{coh}:=\bigcup_{t\in \Delta^*} H^n(X_t,\C),
\end{eqnarray}
by the map which associates the section 
$\left( t\mapsto [\frac{\omega}{\ddd f}|_{X_t}]\in H^n(X_t,\C)\right)$ to the class $[\omega]\in E$, which is represented
by $\omega \in \Omega^{n+1}_{\C^{n+1},0}$.

These sections have moderate growth. We will compare them with
certain distinguished sections of moderate growth:
the elementary sections.

\subsection{Elementary sections}\label{s3.2}

We denoted by $\pi:\www {\Delta^*}\to \Delta^*$ 
the universal covering. 
Now we denote by $H^{flat}$ the space of maps
$\www \Delta^*\to H^{coh}$ 
which are obtained by composing the flat sections of
$\pi^*H^{coh}$ with the canonical projection
$\pi^*H^{coh}\to H^{coh}$. 
This means that $H^{flat}$ is the space of 
``global and multivalued flat sections'' of  $H^{coh}$.
For any $\zeta\in \www{\Delta^*}$, one has a isomorphism
$H^n(X_{\pi(\zeta)},\C)\to H^{flat}$; 
and each structure on $H^n(X_t,\C)$ which is $M$-invariant, 
induces the same structure on $H^{flat}$, 
for example the $\Z$-lattice $H^n(X_t,\Z)$, $M$, $N$, 
$I^{coh}$, $S$, and $h$. 

Starting from the (flat multivalued) sections of $H^{flat}$, 
one can construct holomorphic (univalued) sections 
which have a special behaviour at 0,
the elementary sections: Let $A\in H^{flat}_\lambda$
and let $\alpha\in \Q$  such that $\lambda=e^{-2i\pi\alpha}$;
one defines the elementary section
\begin{eqnarray}\label{3.3}
es(A,\alpha)&=&t^\alpha\exp(\log t \frac{-N}{2i\pi}) A\\
&=& t^\alpha t^{\www N}A \nonumber\\
&=& \sum_{k=0}^n t^\alpha \frac{1}{k!}(\log t)^k{\www N}^k A\nonumber
\end{eqnarray}
where
\begin{eqnarray}\label{3.4}
\www N :=\frac{-N}{2i\pi}.
\end{eqnarray}
One denotes by $C^\alpha$ the space of elementary sections
with fixed $\alpha$. One obtains the decreasing $V$-filtration 
of Malgrange-Kashiwara with  
\begin{eqnarray}\label{3.5}
V^\alpha &:=& \sum_{\beta\geq \alpha}\C\{t\}C^\beta 
=\bigoplus_{\alpha\leq \beta<\alpha+1}\C\{t\}C^\beta \quad \hbox{ and }\\
V^{>\alpha}&:=& \sum_{\beta>\alpha}\C\{t\}C^\beta 
=\bigoplus_{\alpha<\beta\leq \alpha+1}\C\{t\}C^\beta .\nonumber
\end{eqnarray}
These spaces are free $\C\{t\}$-modules of rank $\mu$, and 
$V^{>\alpha}$ for $\alpha\geq -1$ and $V^{\alpha}$ for 
$\alpha>-1$ are also free $\C\{\{\dtu\}\}$-modules of rank 
$\mu$.
They all are subspaces of $\bigcup_{\alpha\in \Q}V^\alpha$,
which is a vector space of dimension $\mu$ over 
$\C\{t\}[t^{-1}]$.
This is the space of germs at 0 of holomorphic sections
of $H^{coh}$ which have moderate growth at 0. 

One has
\begin{eqnarray}\label{3.6}
V^{>-1}\supset E \supset V^{n-1}.
\end{eqnarray}
The first inclusion is a result of B. Malgrange,
the second is an implication of the first one and of the
symmetry of the spetral numbers.
This symmetry has several proofs; 
one of them uses the properties of K. Saito's form
(see \cite{Be}\cite[Remark 10.29 (c)]{He1}), another one
uses the relation between the elementary sections and the
mixed Hodge structure of Steenbrink, which was found
by A.N. Varchenko \cite{AGV}.

\subsection{Barlet's form $\HH_{Barlet}$}\label{s3.3}
D. Barlet \cite{Ba} defined a form 
\begin{eqnarray}\label{3.7}
\www \HH : \Omega^{n+1}_{\C^{n+1},0}\times 
\Omega^{n+1}_{\C^{n+1},0}\to \NN,
\end{eqnarray}
which associates to two $(n+1)$-forms $\omega$ and $\omega'$
the asymptotic expansion
\begin{eqnarray}\label{3.8}
\frac{1}{(2i\pi)^{n}}\int_{X_t} \rho\cdot 
\frac{\omega}{\ddd f}\wedge \oooo{\frac{\omega'}{\ddd f}}
\end{eqnarray}
in 
\begin{eqnarray}\label{3.9}
\NN:=\bigoplus_{\alpha\in\Q,k\in \N} 
\C[[t,\oooo t]]\cdot |t|^{2\alpha}(\log t\oooo t)^k \mod 
\C[[t,\oooo t]];
\end{eqnarray}
here $\rho $ is a $C^\infty$-function which has support 
in a small neighbourhood of 0 and which is equal to 1 
in a smaller neighborhood of 0 in $\C^{n+1}$.

By \cite{Ba}, the form \eqref{3.7} induces a map
\begin{eqnarray}\label{3.10}
\HH_{Barlet} : E\times E \to \NN,
\end{eqnarray}
which is linear in the left argument and semilinear in the
right argument. 
D. Barlet observed that it is related to a hermitian form
on the cohomology of a Milnor fiber by the formula
\eqref{3.14} below. 
The identification of that form with the form $h$ 
was proved by F. Loeser \cite{Loe}
via a compactification of the Milnor fibration.

\subsection{Relation of $\HH_{Barlet}$ with $h$ and $S_Y$}\label{s3.4}
One defines the bundle $H^{coh}_Y$ for a compactification
$f_Y:Y\to \Delta$ of the Milnor fibration analogously to
the formula \eqref{3.2}.
One defines the spaces 
$C^\alpha_Y$, $V^\alpha_Y$ and $V^{>\alpha}_Y$ analogously
to those in subsection \ref{s3.3}.

By \eqref{2.30} one has a canonical projection 
\begin{eqnarray}\label{3.11}
\psi: V^{>-1}_Y \to V^{>-1}
\end{eqnarray}
with kernel $\C\{t\}\cdot \ker (N_Y:C^0_Y\to C^0_Y)$.
One defines $E_Y$ as the preimage of $E$ in $V^{>-1}_Y$.

F. Loeser \cite{Loe} proved that one has the identity
\begin{eqnarray}\label{3.12}
\HH_{Barlet}(a,b) = \frac{1}{(2i\pi)^n}I^{coh}_Y(a_Y,\oooo{b_Y}) 
\quad \mod \C[[t,\oooo t]] 
\end{eqnarray}
for $a,b\in E$ et $a_Y,b_Y\in E_Y$ preimages of $a,b$.
Here the ambiguity of the preimages $a_Y$ and $b_Y$
is compensated by the fact that one considers only the
value modulo $\C[[t,\oooo t]]$.

With  \eqref{2.43}, \eqref{2.44} and the decomposition
of $a$ and $b$ into sums of elementary sections, 
one obtains the formula \eqref{3.14} which relates  
$\HH$ and $h$ and with which D. Barlet \cite{Ba} defined $h$.
One writes $a\in E$ as a sum of elementary sections, 
\begin{eqnarray}\label{3.13}
a&=& \sum_{\alpha>-1} t^\alpha t^{\www N}A(a,\alpha)\\
&=& \sum_{\alpha>-1}\sum_{k=0}^n t^\alpha\frac{1}{k!}(\log t)^k \www N^k
A(a,\alpha), \nonumber
\end{eqnarray}
with $A(a,\alpha)\in H^{flat}_{\exp{(-2i\pi\alpha})}$, 
and in the same way also $b$ with $A(b,\beta)$, and also 
$a_Y$ and $b_Y$ with  
$A_Y(a_Y,\alpha)$ and $A_Y(b_Y,\beta)$. 
Therefore one has 

\begin{eqnarray}
&&\HH_{Barlet}(a,b)\nonumber \\
&=& \sum_{\alpha,\beta>-1} \frac{1}{(2i\pi)^n}I^{coh}_Y
\left(t^\alpha t^{\www N_Y} A_Y(a_Y,\alpha),
\oooo{t^\beta t^{\www N_Y}A_Y(b_Y,\beta)}\right) \nonumber\\
&=& 
\sum_{\alpha,\beta>-1} \sum_{k=0}^n 
t^\alpha\oooo{t^\beta} \frac{1}{k!}(\log t\oooo t)^k 
\frac{1}{(2i\pi)^n}I^{coh}_Y
\left(\www N_Y^k A_Y(a_Y,\alpha),
\oooo{A_Y(b_Y,\beta)}\right)\nonumber\\
&=& 
\sum_{\alpha\notin \N,\beta\notin \N}\sum_{k=0}^n
t^\alpha\oooo{t^\beta} \frac{1}{k!}(\log t\oooo t)^k 
h\left((\www N)^k A(a,\alpha),A(b,\beta)\right) \nonumber\\
&+& 
\sum_{\alpha\in \N,\beta\in \N}\sum_{k=0}^{n-1}
t^\alpha\oooo{t^\beta} \frac{1}{(k+1)!}(\log t\oooo t)^{k+1} \nonumber
h\left((\www N)^k A(a,\alpha),A(b,\beta)\right) \nonumber\\
&=& \sum_{\alpha\notin \N,\beta\notin \N}
h\left(t^\alpha t^{\www N} A(a,\alpha),t^\beta t^{\www N}A(b,\beta)\right)
+\sum_{\alpha\in \N,\beta\in \N}\sum_{k=0}^{n-1} ...
\label{3.14}
\end{eqnarray}

\subsection{K. Saito's bilinear form $K_{Saito}$ and the 
bilinear form $J_{Grothendieck}$}\label{s3.5}
The Grothendieck residue is a linear form on the Jacobi 
algebra, defined by 

\begin{eqnarray} \label{3.15}
&&\OO_{\C^{n+1},0}/
\left(\frac{\paa f}{\paa x_0},...,\frac{\paa f}{\paa x_n}\right) \to \C\\
&&g\mapsto \Res_0\left[\frac{g\ddd x}
{\frac{\paa f}{\paa x_0}\cdot ...\cdot \frac{\paa f}{\paa x_n}}\right]
:=\frac{1}{(2i\pi)^{n+1}}\int_{\Gamma_\varepsilon} 
\frac{g\ddd x}{\frac{\paa f}{\paa x_0}\cdot ...\cdot \frac{\paa f}{\paa x_n}},
\nonumber
\end{eqnarray}
where 
$\Gamma_\varepsilon :=\{ x\in \C^{n+1}\ |\ 
|\frac{\paa f}{\paa x_i}\ |=\varepsilon\quad \forall\ i\}$.
It does not depend on $\varepsilon$, but it depends on the
coordinates $x_0,...,x_n$.

K. Saito remarked that it induces a bilinear form $J_f$
on the $\mu$ dimensional vector space
(which is a free module of rank 1 of the Jacobi algebra)
\begin{eqnarray}\nonumber
\Omega_f&:=& \Omega_{\C^{n+1},0}^{n+1}/\ddd f\wedge \Omega_{\C^{n+1},0}^n\\
&=& E/\dtu E \label{3.16},
\end{eqnarray}
which is defined by 
\begin{eqnarray}\label{3.17}
J_f:\Omega_f\times\Omega_f &\to& \C,\\
(g_1\ddd x,g_2\ddd x) &\mapsto& \Res_0\left[\frac{g_1g_2\ddd x}
{\frac{\paa f}{\paa x_0}\cdot ...\cdot \frac{\paa f}{\paa x_n}}\right] .
\nonumber
\end{eqnarray}
Furthermore, it is independent of the coordinates 
$x_0,...,x_n$.
The form $J_f$ is symmetric and nondegenerate;
the nondegeneracy is nontrivial and was shown first by
A. Grothendieck.

Motivated by \eqref{3.17} and \eqref{3.16}, 
K. Saito was able to define a generalization of $J_f$ 
(see \cite{SaK}\cite{Na}\cite{Be}). He constructed a form
\begin{eqnarray}\label{3.18}
K_f: E\times E\to \dt^{-n-1}\C\{\{\dtu\}\}
\end{eqnarray}
with the following four properties:
\begin{eqnarray}\label{3.19}
&i)& \quad K_f(\dtu a,b)=\dtu K_f(a,b) = K_f(a,-\dtu b);\hspace*{1cm}\\
&ii)& \quad K_f(t\, a,b)-K_f(a,t\, b) = [t,K_f(a,b)];\hspace*{1cm}\label{3.20}
\end{eqnarray}
where $[t,\dt^{-k}]= k\dt^{-k-1}$;
one defines $K_f^{(-k)}:E\times E\to \C$ for $k\in \N$ by
\begin{eqnarray}\label{3.21}
K_f =\sum_{k=0}^\infty K_f^{(-k)}\cdot \dt^{-n-1-k}.
\end{eqnarray} 
$iii)$ \quad $K_f^{(-k)}$ is $(-1)^k$ symmetric; \\
$iv)$ \quad the form $K_f^{(0)}$ satisfies
\begin{eqnarray}\label{3.22}
K_f^{(0)}(E,\dtu E)=K_f^{(0)}(\dtu E,E)=0
\end{eqnarray}
and it induces the form $J_f$ on the quotient 
$\Omega_f=E/\dtu E$.

Here we will not discuss the (algebraic) definition of
K. Saito, but we will give below an analytic formula
due to A.N. Varchenko.

\subsection{Relation of $K_{Saito}$ with $S$ and $S_Y$}\label{s3.6}
A.N. Varchenko \cite{Va} considered like J.H.M. Steenbrink and
J. Scherk and later F. Loeser \cite{Loe} (cf. \ref{s3.4})
a compactification $f_Y:Y\to \Delta$ of the Milnor fibration.
He defined a series of pairings 
\begin{eqnarray}\label{3.23}
\www K_f^{(-k)}:V^{>-1}\times V^{>-1}\to \C\quad \hbox{ for }
k\in \Z_{\geq -n}.
\end{eqnarray}
With the notations of \ref{s3.4} they are given by 
\begin{eqnarray}\label{3.24}
\www K_f^{(-k)}(a,b) := (-1)^{n(n-1)/2}\frac{1}{(2i\pi)^{n}}
\res_0 I^{coh}_Y(\nnn_{\dt}^{n+k} a_Y,b_Y)
\end{eqnarray}
(here we calculated certain constants which had not been
calculated in \cite{Va}).
The constants are chosen such that $\www K_f^{(0)}$ 
restricted to $E\subset V^{>-1}$ induces $J_f$ on $\Omega_f$. 
One can prove
\begin{eqnarray}\label{3.25}
\www K_f^{(-k)}|_{E\times E} = K_f^{(-k)} 
\quad\hbox{ for }k\in \N.
\end{eqnarray}
The proof uses several arguments:
\begin{list}{}{}
\item[a)] 
The form $K_f$ is the restriction of a form $K_F$
for the Gauss-Manin system of a semiuniversal unfolding $F$;
the form $K_F$ was also defined by K. Saito.
\item[b)] 
The form $K_F$ can be characterized uniquely if one 
generalizes the properties above of $K_f$.
This results from the fact that the Gauss-Manin system
of $F$ is simple holonomic as a microdifferential system
(see \cite[(2.7.11)]{SaM}).
\item[c)]
The form $\www K_f^{(-k)}$ for $k\in\N$ can also be extended
to the Gauss-Manin system of $F$, 
keeping the same properties
\end{list}
Thus thanks to \eqref{3.25}, one can also define 
$K_f^{(-k)}$ for $k\in \{-n,...,-1\}$, 
and all these forms are defined on $V^{>-1}$;
one obtains therefore a form 
\begin{eqnarray}\label{3.26}
K_f:V^{>-1}\times V^{>-1}\to \dtu\C\{\{\dtu\}\}.
\end{eqnarray}
Thanks to \eqref{3.24}, one sees that for $\alpha>-1$ and 
$\beta>-1$
\begin{eqnarray}\label{3.27}
&&K_f:C^\alpha\times C^\beta\to 0 \quad\hbox{ if }\alpha+\beta\notin\Z\\
&&K_f:C^\alpha\times C^\beta\to \C\cdot\dt^{-\alpha-\beta-2}
\quad\hbox{ if }\alpha+\beta\in\Z \label{3.28}
\end{eqnarray}
and that the form in \eqref{3.28} is nondegenerate and
$(-1)^{\alpha+\beta+n+1}$-symmetric.
In fact, starting from \eqref{3.24} and \ref{s2.8}, one can
prove the following formulas, with 
$A\in H^{flat}_{\exp{(-2i\pi\alpha)}}$,
$B\in H^{flat}_{\exp{(-2i\pi\beta)}}$ and the notations
of \ref{s3.2}
\begin{eqnarray}\label{3.29}
K_f(es(A,\alpha),es(A,\beta))&=&\frac{1}{(2i\pi)^n}S(A,B)\cdot \dt^{-1}\\ 
  && \hbox{ if } \alpha,\beta\in ]-1,0[,\ \alpha+\beta=-1,\nonumber \\
K_f(es(A,\alpha),es(A,\beta))&=&\frac{-1}{(2i\pi)^{n+1}}S(A,B)\cdot\dt^{-2}\label{3.30}\\
  && \hbox{ if } \alpha=\beta=0.\nonumber
\end{eqnarray}
The formulas \eqref{3.29}, \eqref{3.30} and  \eqref{3.19}
give a topological definition of $K_f$.\\
The form $K_f^{(-1)}$ has interesting properties:
it is antisymmetric, it induces a nondegenerate form
on the quotient $V^{>-1}/V^{n-1}$, 
and the Brieskorn lattice $E$ induces a subspace $E/V^{n-1}$ 
of $V^{>-1}/V^{n-1}$ which is Lagrangian with respect 
to $K_f^{(-1)}$.

\subsection{Fourier-Laplace transformation}\label{s3.7}
Also for the cohomology bundle $H^{Lef}_\C$ 
of the Lefschetz thimbles (see \ref{s2.6}), one can define
the notion of an elementary section, as well as the spaces
$C^\alpha_{Lef}$ and a space $H^{flat}_{Lef}$, 
in the same way as it was done for the bundle $H^{coh}$ in \ref{s3.2}. 

The isomorphism $\paa_*$ in \eqref{2.22}
induces a topological isomorphism 
\begin{eqnarray}\label{3.31}
\paa^*_{Lef}:H^{flat}\to H^{flat}_{Lef}.
\end{eqnarray}
One has a Fourier-Laplace transformation $FL$
\begin{eqnarray}\label{3.32}
FL: V^{mod}&\stackrel{\cong}{\longrightarrow}&V^{mod}_{Lef},\\
C^\alpha&\stackrel{\cong}{\longrightarrow}&C^{\alpha+1}_{Lef}
\quad\hbox{ pour }\alpha>-1,\nonumber
\end{eqnarray}
between the spaces 
\begin{eqnarray}\label{3.33}
V^{mod}&:=&\bigoplus_{-1<\alpha\leq 0}\C[t]C^\alpha \quad 
\hbox{ and }\\
V^{mod}_{Lef}&:=&\bigoplus_{0<\alpha\leq 1}\C[z]C^\alpha_{Lef}.
\label{3.34}
\end{eqnarray}
One can define $FL$ either algebraically by 
\begin{eqnarray}\label{3.35}
\dtu &\mapsto& z,\\
\dt  &\mapsto& z^{-1},\nonumber \\
t    &\mapsto& -\paa_{z^{-1}}=z^2\paa_z,\nonumber
\end{eqnarray}
(see \cite{Sab2}) or analytically by the following 
formula, for $a(t)\in V^{mod}$
\begin{eqnarray}\label{3.36}
FL(a(t))(z):= \int_0^{\infty\cdot \frac{z}{|z|}}e^{-\frac{t}{z}}\cdot 
\paa^*_{Lef} a(t) \ddd t.
\end{eqnarray}
The Brieskorn lattice $E$ gives via $FL$ a subspace 
\begin{eqnarray}\label{3.37}
G_0:= FL(E\cap V^{mod}) \subset V^{mod}_{Lef},
\end{eqnarray}
which is also called Brieskorn lattice in \cite{Sab1}.
It gives an extension of the bundle $H^{Lef}_\C\to \C^*$ 
to a bundle on $\C$.

\subsection{Relation of $K_{Saito}$ with $P_{Pham}$}\label{s3.8}
On the cohomology bundle $H^{Lef}_\C$ of the Lefschetz 
thimbles, one has the pairing $P_{Pham}$ of Pham (\eqref{2.29}).
Pham \cite{Ph2} proved the formula 
\begin{eqnarray}\label{3.38}
K_f(a,b) \ ''='' P_{Pham}(FL(a),FL(b))\quad 
\hbox{ for }a,b\in E\cap V^{mod}.
\end{eqnarray}
Here $''=''$ means that one identifies a polynomial in $\dtu$ 
on the left hand side with a polynomial in $z$ on the 
right hand side.

Pham's proof uses the same arguments a) - c)  for \ref{s3.6}
as the proof of \eqref{3.25}: The form $P_{Pham}$ can be 
defined for a semiuniversal unfolding $F$ of $f$
without difficulties, and after the identification $''=''$,
it satisfies the same properties as $K_F$. 

But one could replace this proof by the characterization
of $K_f$ with \eqref{3.29}, \eqref{3.30} and \eqref{3.19},
by a comparison of the topological forms $P_{Pham}$ and $S$
and by explicit calculations which yield the following formulas
(this is explained without details in \cite[7.2]{He2},
all details are given in \cite[ch. 4 and 5]{BH}):
\begin{eqnarray}\label{3.39}
&&FL(es(A,\alpha-1))(z)=es(G^{(\alpha)}A,\alpha)(z)
\end{eqnarray}
for $\alpha>0, A\in H^{flat}_{\exp(-2i\pi\alpha)}$, where 
\begin{eqnarray}\label{3.40}
G^{(\alpha)}&=& \paa_{Lef}^*\circ ''\Gamma(\alpha\cdot \id+\www N)''\\
&=& \paa_{Lef}^*\circ \sum_{k\geq 0}\frac{1}{k!}\Gamma^{(k)}(\alpha)
\www N^k : H^{flat}\to H^{flat}_{Lef};\nonumber
\end{eqnarray}
and
\begin{eqnarray}\label{3.41}
P_{Pham}\left(es(G^{(\alpha)}A,\alpha),es(G^{(\beta)}B,\beta)\right) 
= z\frac{1}{(2i\pi)^n}\cdot S(A,B)
\end{eqnarray}
for $\alpha,\beta\in ]0,1[$ with $\alpha+\beta=1$, 
$A\in H^{flat}_{\exp(-2i\pi\alpha)}$, $B\in H^{flat}_{\exp(-2i\pi\beta)}$ 
and
\begin{eqnarray}\label{3.42}
P_{Pham}\left(es(G^{(1)}A,1),es(G^{(1)}B,1)\right) 
= z^2\frac{-1}{(2i\pi)^{n+1}}\cdot S(A,B)
\end{eqnarray}
for $\alpha=\beta=1$ et 
$A\in H^{flat}_1$, $B\in H^{flat}_1$.

This proof shows that \eqref{3.38} is true for $a,b\in V^{mod}$.

\end{document}